\documentstyle[amssymb,12pt,leqno,]{article}

\textwidth=16.50 true cm
\textheight=24.00 true cm

\hoffset=-1.5cm
\voffset=-2cm
\parindent=25pt




\newcounter{draw}[section]

\renewcommand\thedraw{\thesection.\the\nf}

\global\newcount \nf



\newcommand{\bg}{\begin}
\newcommand{\ed}{\end}

\newcommand{\dy}{\displaystyle}

\newcommand{\tx}{{\cal X}\kern -0.30cm\raise 0.04cm
\hbox{$-$}}
\newcommand{\bgpf}[1]{\noindent {\bf Proof{#1}:}}


\headsep=15pt

\title{Lower bounds for index of Wente tori}
\author{Levi Lopes de Lima, Vicente Francisco de Sousa Neto and 
Wayne Rossman}
\date{}

\bg{document}
\maketitle
\baselineskip=14pt
\abovedisplayskip=1pc
\belowdisplayskip=1cm

\begin{center}
Dedicated to Katsuhiro Shiohama on the occasion of his sixtieth 
birthday.  
\end{center}

\begin{abstract}
We show numerically that any of the constant mean curvature tori first 
found by Wente must have index at least eight.
\end{abstract}

\section{Introduction}
The Hopf conjecture asked if all closed surfaces 
immersed in ${\Bbb R}^3$ with constant mean curvature $H$ must 
be round spheres.  It was proven true when either the surface has 
genus zero by Hopf himself \cite{H}, or the immersion is 
actually an embedding by Alexandrov \cite{H}.  However, it 
does not hold in general, and the first counterexamples, of genus $1$, 
were found by Wente \cite{We}.  Abresch \cite{A} and Walter \cite{Wa} made 
more explicit 
descriptions for these surfaces of Wente, which all have one family of 
planar curvature lines \cite{Sp}.  We call these surfaces the 
{\em original Wente tori}.

Constant mean curvature 
surfaces are critical for area, but not necessarily area minimizing, for 
all compactly supported volume-preserving variations.  
Hence the index -- loosely speaking, the {\em dimension} of area-reducing 
volume-preserving variations, to be 
defined in Section 3 -- can be positive.  If it is zero, the 
surface is stable.  Do Carmo and Peng \cite{CP} showed that the only 
complete stable minimal surface is a plane.  
Fischer-Colbrie \cite{FC} showed that a complete minimal 
surface in ${\Bbb R}^3$ has 
finite index if and only if it has finite total curvature, and 
that the catenoid and Enneper's surface have index $1$.  Likewise, for 
surfaces with constant mean curvature 
$H \neq 0$, Barbosa and Do Carmo \cite{BC} showed 
that only round spheres are stable, and 
Lopez and Ros \cite{LR} and Silveira \cite{Si} independently showed that 
they have finite index if and only if they are compact.  This leaves open 
the question of whether there exist surfaces with constant mean curvature 
$H \neq 0$ and low positive index, for example with index $1$.  

The third author \cite{R1}, \cite{R2} showed numerically that 
the most natural candidates for unstable surfaces of constant mean 
curvature $H \neq 0$ with low index -- the original 
Wente tori -- all have index at least $7$, 
and with a numerical experiment suggested that the sharpest lower bound 
is either $8$, $9$, or $10$, and is most likely $9$.  
This leads one to conjecture that all closed surface with 
constant mean curvature $H \neq 0$ have index at least $9$.  

The purpose of 
this article is to show that the original Wente tori all have index 
at least $8$, improving the lower bound of \cite{R1}, \cite{R2}.  
The final part of our argument relies on numerics.  

\section{The original Wente tori}

In this section we shall give a brief description of the original Wente 
tori, based on \cite{Wa}.  Later, we shall assume that the mean curvature 
$H$ is $1/2$, but in this and the next section we shall only assume that 
$H$ is a nonzero constant.  

Let $\tx: {\Bbb C}/\Gamma 
\longrightarrow{\Bbb R}^3$ be a conformal immersion of class
$C^\infty$ where ${\Bbb C}/\Gamma$ is a compact $2$-dimensional torus 
determined by the $2$-dimensional lattice $\Gamma$.  
Note that $(x,y)$ are then isothermal coordinates on ${\Bbb C}/\Gamma$.  
The fundamental forms and the Gauss and mean curvature functions are 
\begin{center}
$I=E(dx^2+dy^2)$, $II=Ldx^2+2Mdx\,dy+Ndy^2$, 
$K=\frac{LN-M^2}{E^2}$, $H=\frac{L+N}{2E}$. 
\end{center}
Since $H$ is constant, the Hopf differential $\Phi\,dz^2$ is holomorphic, 
where $\Phi=\frac{1}{2}(L-N)-iM$ and 
$z=x+iy$.  Thus $\Phi$ is constant and 
$\tx$ has no umbilics points.  Moreover, by a change of the 
coordinates $(x,y)$, we may assume $\Phi=1$ and so $M=0$, 
$L=e^F+1$, $N=e^F-1$, and $(x,y)$ become curvature line parameters, 
where $F: {\Bbb C}/\Gamma \longrightarrow{\Bbb R}$ is defined by $HE=e^F$. 
We have the equations of Gauss and Weingarten:
\bg{equation}
\bg{array}{l}
\dy\tx_{xx}=\frac{1}{2}F_x\tx_x-\frac{1}{2}F_y\tx_y-(e^F+1){\cal N} \; , 
\; \; \; 
\dy\tx_{yy}=-\frac{1}{2}F_x\tx_x+\frac{1}{2}F_y\tx_y-(e^F-1){\cal N} \; , \\
\dy\tx_{xy}=\frac{1}{2}F_y\tx_x+\frac{1}{2}F_x\tx_y \; , \; \; \; 
\dy {\cal N}_{x}=H(1+e^{-F})\tx_x \; , \;\;\; 
\dy {\cal N}_{y}=H(1-e^{-F})\tx_y
\ed{array}
\ed{equation}
\bg{equation}
\Delta F+4H\,\sinh\,F=0 \; , 
\ed{equation}
where $\Delta=\frac{\partial^2}{\partial
x^2}+\frac{\partial^2}{\partial y^2}$ and ${\cal N}: {\Bbb C}/\Gamma 
\longrightarrow R^3$ is the unit normal vector field, i.e. the Gauss map.  
Therefore the problem of finding constant mean curvature immersed tori in
${\Bbb R}^3$ reduces to solving the PDE system (1) and (2) by 
real analytic functions $F,{\cal N},\tx$ defined on ${\Bbb R}^2$ and doubly
periodic with respect to some fundamental lattice $\Gamma\subset
{\Bbb R}^{2}$. 

In the case of the original Wente tori, in Walter's notation, the 
solution $F$ of (2) is:
\bg{equation}
\tanh\left(\frac{F}{4}\right)=\gamma\cdot\overline{\gamma}\cdot
cn_k \left( \alpha x \right) cn_{\overline{k}} \left( \overline{\alpha}y 
\right) \; , 
\ed{equation}
where $cn_k$ denotes the Jacobi amplitudinus cosinus 
function with modulus $k$, and $k=\sin\,\theta,\
\overline{k}=\sin\,\overline{\theta}$, for
$\theta,\overline{\theta}\in(0,\pi/2)$ and
$\theta+\overline{\theta}<\pi/2$, and 
\begin{center}
$\gamma=\sqrt{\tan\,\theta}$, 
$\overline{\gamma}=\sqrt{\tan\,\overline{\theta}}$, 
$\alpha=\sqrt{4H\frac{\sin\,2\overline{\theta}}
{\sin\,2(\theta+\overline{\theta})}}$, 
$\overline{\alpha}=\sqrt{4H\frac{\sin\,2\theta}
{\sin\,2(\theta+\overline{\theta})}} \;\;\;\;\; $.
\end{center}

\bg{lema} {\em (\cite{A}, \cite{Wa})}.  
The set of all original Wente tori are in a one-to-one
correspondence with the set of reduced fractions $\ell/n\in(1,2)$. 
\ed{lema}


For each $\ell/n$, we call the corresponding 
Wente torus ${\cal W}_{\ell/n}$.  Following Walter's notation, each 
${\cal W}_{\ell/n}$ has either one or two planar geodesic loops in the 
central symmetry plane: two loops 
if $\ell$ is odd, and one loop if $\ell$ is even.  Each 
loop can be partitioned into $2 n$ congruent curve segments, 
and $\ell$ is the total winding order of the Gauss map along each 
loop.  

The conditions for double periodicity of the position vector function
$\tx$ are expressed in terms of $\theta$ and
$\overline{\theta}$. Walter determined
that there is exactly one 
\[ \overline{\theta} \cong 65.354955^\circ \] 
that solves one period problem.  The other period problem is solved with the 
correct choice of $\theta \in (0,\ (\pi/2)-\overline{\theta})$, and, for any 
$\ell/n\in (1,2)$, this correct choice is the unique solution $\theta$ of
\bg{equation}
\int^{\pi/2}_0\frac{1+\tan \theta\ \tan 
\overline{\theta}\ \cos^2\varphi}
{1-\tan \theta\ \tan \overline{\theta}\ \cos^2\varphi}
\frac{d\varphi}{\sqrt{1-\sin^2\theta\ \sin^2\varphi}}=
\frac{\ell}{n}\frac{\pi}{2}
\sqrt{\frac{\sin\,2\overline{\theta}}{\sin\,2(\theta+\overline{\theta})}} \; . 
\ed{equation}
For any $\ell/n\not\in(1,2)$, there is no solution $\theta\in 
(0,(\pi/2)-\overline{\theta})$ of (4).  
In Table 2 we give some values of $\theta$ with respect to $\ell/n$.

Now, if $x_{\ell n}$ (resp. $y_{\ell n}$) denotes the length of the period 
of $cn_k(\alpha x)$ (resp. $cn_{\overline{k}}(\overline{\alpha}y)$), then
we have the following lemma:  

\bg{lema} {\em (\cite{Wa})}.  $\tx:{\Bbb C}/\Gamma\longrightarrow 
W_{\ell/n} \subset {\Bbb R}^3$ is a 
conformal immersion {\em (}$W_{\ell/n}$ denotes the image of $\tx${\em )}, 
where $$\Gamma=\mbox{span}_{\Bbb Z}\{(nx_{\ell n},0),(0,y_{\ell n})\} 
\quad\mbox{when}\; \ell \; \mbox{is odd, and}$$
$$\Gamma=\mbox{span}_{\Bbb
Z}\left\{\left(\frac{nx_{\ell n}}{2},\frac{y_{\ell
n}}{2}\right),(0,y_{\ell n})\right\}
\quad\mbox{when} \; \ell \; \mbox{is even}.$$
The curves $\{[x_0,y]|x_0=\mbox{constant}\}$ are mapped by $\tx$ to
planar curvature lines of $W_{\ell/n}$. 
\ed{lema}

The lengths $x_{\ell n}$ and $y_{\ell n}$ can be computed as follows: 
\bg{equation} 
x_{\ell n}=\frac{4}{\alpha}\int^{\pi/2}_0 
\frac{d\varphi}{\sqrt{1-k^2\sin^2\varphi}} \; , \; \; \; 
 y_{\ell n}=\frac{4}
{\overline{\alpha}}\int^{\pi/2}_0
\frac{d\varphi}{\sqrt{1-{\overline{k}}^2\sin^2\varphi}} \; . 
\ed{equation}

\section{The definition of index and preliminary results} 

The Jacobi operator associated to ${\cal W}_{\ell/n}$ is 
$-\Delta_I-|II|^2$ on ${\Bbb C}/\Gamma$, with $|II|^2 = 
E^{-2} (L^2+2 M^2+N^2)=2H^2 (1+e^{-2F})$ and $\Delta_I$ the 
Laplace-Beltrami operator associated to the metric $I$.  
The corresponding quadratic form is 
\bg{equation}
Q(u,u)=\int_{{\Bbb C}/\Gamma}u{\cal L}(u) \, dxdy,
\ed{equation}
where \[ {\cal L}u=-\Delta u-Vu \;\;\;\;\;\; \mbox{with} \;\;\;\;\;\; 
V=4H\cosh(F) \] and $\Delta$ the Euclidean Laplacian.  
Note that in equation (6), we are integrating with respect to the
flat metric on ${\Bbb C}/\Gamma$. 

Consider a smooth volume-preserving variation $\tx_t$ of the immersion 
$\tx$ with parameter 
$t$ so that $\tx_0$ is the surface ${\cal W}_{\ell/n}$.  By reparametrizing 
the surfaces of the variation, we may assume 
that the variation vector field at $t=0$ is $u {\cal N}$ for some 
$u \in C^\infty({\Bbb C}/\Gamma)$.  Then 
\begin{center}
$\left. 
\frac{\partial}{\partial t} \mbox{area}(\tx_t) 
\right|_{t=0}=0$ and 
$\left. 
\frac{\partial^2}{\partial t^2} \mbox{area}(\tx_t)\right|_{t=0}=
Q(u,u)$.  
\end{center}
Furthermore, the 
volume-preserving condition implies $\int_{{\Bbb C}/\Gamma} u\,dA=0$. 
Thus, if 
\begin{center}
${\cal V}=\left\{u\in C^{\infty}({\Bbb C}/\Gamma); \int_{{\Bbb C}/\Gamma} 
u\,dA=0\right\}$, 
\end{center}
then we can give the following definition (see \cite{BC}):

\bg{dft}
{\em 
We define $\mbox{Ind}(\tx ({\Bbb C}/\Gamma))$, the 
{\it index} of the immersion 
$\tx$ of ${\Bbb C}/\Gamma$, to be the maximum of the 
dimensions of the subspaces of ${\cal V}$ restricted to which $Q$ is
negative definite. 
}
\ed{dft}

Since the first derivative of area is zero, and the second derivative 
is $Q(u,u)$, the index in a sense measures the amount of area-reducing 
volume-preserving variations.  

Let $L^2=L^2({\Bbb C}/\Gamma)=\{u\in 
C^\infty({\Bbb C}/\Gamma)|\int_{{\Bbb C}/\Gamma}u^2dx\,dy<\infty\}$
provided with the inner product $\langle u,v\rangle_{L^2}=\int_{{\Bbb
C}/\Gamma}uv\,dx\,dy$. It follows from the standard spectral
theorem that the operator ${\cal L}=-\Delta-V$ on ${\Bbb C}/\Gamma$ has a 
discrete spectrum of eigenvalues
\begin{center}
$\beta_1 < \beta_2\leq\cdots\nearrow\infty$ \end{center} 
and corresponding eigenfunctions 
\begin{center}
$\nu_1,\nu_2\cdots\in C^\infty({\Bbb C}/\Gamma)$, 
\end{center}
which form an orthonormal basis for $L^2$. Moreover, we have the
following variational characterization for the eigenvalues:
\begin{center}
$\beta_j=\inf_{V_j}\left(\sup_{\phi\in V_j,\ ||\phi||_{L^2}=1}
\int_{{\Bbb C}/\Gamma}\phi{\cal L}\phi dx\,dy\right)$, \end{center}
where $V_j$ runs through all $j$ dimensional subspaces of
$C^\infty({\Bbb C}/\Gamma)$. 

\bg{lema} {\em (\cite{R1}, \cite{R2})}.  If ${\cal L}$ has $k$ negative 
eigenvalues, then {\em Ind}$(W_{\ell/n})$ is either $k$ or $k-1$.  
Furthermore, if there exists a 
subspace ${\mathcal{S}} \subset L^2$ such that ${\mathcal{S}} \subset 
C^\infty({\Bbb C}/\Gamma)$  and {\em dim}$({\mathcal{S}})=k$ and 
$Q$ restricted to 
$\mathcal{S}$ is negative definite, then {\em Ind}$(W_{\ell/n})\geq k-1$. 
\ed{lema}

\medskip

%
%

By Lemma 3, our goal becomes to compute the number of negative
eigenvalues of ${\cal L}$.

Now, we use a convenient fact: For the flat torus  ${\Bbb C}/\Gamma$, with
$\Gamma=\mbox{span}_{\Bbb Z}\{(a_1,a_2),(b_1,b_2)\}$, the complete set
of eigenvalues of $-\Delta u_i=\alpha_iu_i$ are 
\begin{center}
$\frac{4\pi^2}{(a_1b_2-a_2b_1)^2}((m_2b_2-m_1a_2)^2+(m_1a_1-m_2b_1)^2)$, 
\end{center}
with corresponding orthonormal eigenfunctions
\begin{center}
$c_{m_1,m_2}\cdot(\sin\ \mbox{or}\ \cos)
\left(\frac{2\pi}{a_1b_2-a_2b_1}((m_2b_2-m_1a_2)x+
(m_1a_1-m_2b_1)y)\right)$, \end{center}
where $m_1,m_2\in{\Bbb Z}$, 
$c_{m_1,m_2}=\sqrt{2/(a_1b_2-a_2b_1)}$ if 
$|m_1|+|m_2|>0$, $c_{0,0}=\sqrt{1/(a_1b_2-a_2b_1)}$.  
With the aid of Lemma 2 we list $17$ of the $\alpha_i$ and $u_i$ in 
Table 1.  

With the orderings for the eigenvalues as chosen in 
Table 1, we do not necessarily have
$\alpha_i\leq\alpha_j$ for $i\leq j$. However, we still have
$\alpha_i\nearrow\infty$ as $i\nearrow\infty$. Choose
$\alpha_{\rho_{\ell/n}(1)},\alpha_{\rho_{\ell/n}(2)},\cdots$ the complete 
set of eigenvalues with multiplicity 1 of the operator $-\Delta$ on the flat
torus ${\Bbb C}/\Gamma$ reordered by the permutation $\rho_{\ell/n}$ of 
$\Bbb N$ so that $\alpha_{\rho_{\ell/n}(1)} < 
\alpha_{\rho_{\ell/n}(2)}\leq\cdots\nearrow\infty$.

\bg{tabela}{
\bg{tabular}{c|c|c|c}
eigenvalues & eigenfunctions & eigenvalues & eigenfunctions\\ 
$\alpha_i$ for $\ell$ odd & $u_i$ for $\ell$ odd &
$\alpha_i$ for $\ell$ even & $u_i$ for $\ell$ even \\ \hline
$\alpha_1=0$ & $u_1=\frac{1}{\sqrt{nx_{\ell n}y_{\ell n}}}$ &
$\alpha_1=0$ & $u_1=\frac{\sqrt{2}}{\sqrt{nx_{\ell n}y_{\ell n}}}$ \\ \hline
$\alpha_2=\frac{4\pi^2}{n^2x^2_{\ell n}}$ & $u_2=\frac{\sin(2\pi
x/(nx_{\ell n}))}{\sqrt{nx_{\ell n}y_{\ell n}/2}}$ &
$\alpha_2=\frac{16\pi^2}{n^2x^2_{\ell n}}$ & $u_2=\frac{\sin(4\pi
x/(nx_{\ell n}))}{\sqrt{nx_{\ell n}y_{\ell n}/4}}$ \\ \hline
$\alpha_3=\frac{4\pi^2}{n^2x^2_{\ell n}}$ & $u_3=\frac{\cos(2\pi
x/(nx_{\ell n}))}{\sqrt{nx_{\ell n}y_{\ell n}/2}}$ &
$\alpha_3=\frac{16\pi^2}{n^2x^2_{\ell n}}$ & $u_3=\frac{\cos(4\pi
x/(nx_{\ell n}))}{\sqrt{nx_{\ell n}y_{\ell n}/4}}$ \\ \hline
$\alpha_4=\frac{4\pi^2}{y^2_{\ell n}}$ & $u_4=\frac{\sin(2\pi
y/y_{\ell n})}{\sqrt{nx_{\ell n}y_{\ell n}/2}}$ &
$\alpha_4=\frac{4\pi^2}{n^2x^2_{\ell n}}+\frac{4\pi^2}{y^2_{\ell n}}$ & 
$u_4=\frac{\sin(\frac{2\pi x}{nx_{\ell n}}+\frac{2\pi y}{y_{\ell n}})}
{\sqrt{nx_{\ell n}y_{\ell n}/4}}$ \\ \hline
$\alpha_5=\frac{4\pi^2}{y^2_{\ell n}}$ & $u_5=\frac{\cos(2\pi
y/y_{\ell n})}{\sqrt{nx_{\ell n}y_{\ell n}/2}}$ &
$\alpha_5=\frac{4\pi^2}{n^2x^2_{\ell n}}+\frac{4\pi^2}{y^2_{\ell n}}$ & 
$u_5=\frac{\cos(\frac{2\pi x}{nx_{\ell n}}+\frac{2\pi y}{y_{\ell n}})}
{\sqrt{nx_{\ell n}y_{\ell n}/4}}$ \\ \hline
$\alpha_6=\frac{16\pi^2}{n^2x^2_{\ell n}}$ & $u_6=\frac{\sin(4\pi
x/(nx_{\ell n}))}{\sqrt{nx_{\ell n}y_{\ell n}/2}}$ &
$\alpha_6=\frac{4\pi^2}{n^2x^2_{\ell n}}+\frac{4\pi^2}{y^2_{\ell n}}$ & 
$u_6=\frac{\sin(\frac{2\pi x}{nx_{\ell n}}-\frac{2\pi y}{y_{\ell n}})}
{\sqrt{nx_{\ell n}y_{\ell n}/4}}$ \\ \hline
$\alpha_7=\frac{16\pi^2}{n^2x^2_{\ell n}}$ & $u_7=\frac{\cos(4\pi
x/(nx_{\ell n}))}{\sqrt{nx_{\ell n}y_{\ell n}/2}}$ &
$\alpha_7=\frac{4\pi^2}{n^2x^2_{\ell n}}+\frac{4\pi^2}{y^2_{\ell n}}$ & 
$u_7=\frac{\cos(\frac{2\pi x}{nx_{\ell n}}-\frac{2\pi y}{y_{\ell n}})}
{\sqrt{nx_{\ell n}y_{\ell n}/4}}$ \\ \hline
$\alpha_8=\frac{4\pi^2}{n^2x^2_{\ell n}}+\frac{4\pi^2}{y^2_{\ell n}}$ & 
$u_8=\frac{\sin(\frac{2\pi x}{nx_{\ell n}}+\frac{2\pi y}{y_{\ell n}})}
{\sqrt{nx_{\ell n}y_{\ell n}/2}}$ &
$\alpha_8=\frac{16\pi^2}{y^2_{\ell n}}$ & 
$u_8=\frac{\sin(4\pi y/y_{\ell n})}{\sqrt{nx_{\ell n}y_{\ell n}/4}}$ \\
\hline 
$\alpha_9=\frac{4\pi^2}{n^2x^2_{\ell n}}+\frac{4\pi^2}{y^2_{\ell n}}$ & 
$u_9=\frac{\cos(\frac{2\pi x}{nx_{\ell n}}+\frac{2\pi y}{y_{\ell n}})}
{\sqrt{nx_{\ell n}y_{\ell n}/2}}$ &
$\alpha_9=\frac{16\pi^2}{y^2_{\ell n}}$ & 
$u_9=\frac{\cos(4\pi y/y_{\ell n})}{\sqrt{nx_{\ell n}y_{\ell n}/4}}$ \\
\hline 
$\alpha_{10}=\frac{4\pi^2}{n^2x^2_{\ell n}}+\frac{4\pi^2}{y^2_{\ell n}}$ & 
$u_{10}=\frac{\cos(\frac{2\pi x}{nx_{\ell n}}-\frac{2\pi y}{y_{\ell n}})}
{\sqrt{nx_{\ell n}y_{\ell n}/2}}$ &
$\alpha_{10}=\frac{64\pi^2}{n^2x^2_{\ell n}}$ & 
$u_{10}=\frac{\sin(8\pi x/(nx_{\ell n}))}{\sqrt{nx_{\ell n}y_{\ell n}/4}}$ \\ \hline
$\alpha_{11}=\frac{4\pi^2}{n^2x^2_{\ell n}}+\frac{4\pi^2}{y^2_{\ell n}}$ & 
$u_{11}=\frac{\cos(\frac{2\pi x}{nx_{\ell n}}-\frac{2\pi y}{y_{\ell n}})}
{\sqrt{nx_{\ell n}y_{\ell n}/2}}$ &
$\alpha_{11}=\frac{64\pi^2}{n^2x^2_{\ell n}}$ & 
$u_{11}=\frac{\cos(8\pi x/(nx_{\ell n}))}{\sqrt{nx_{\ell n}y_{\ell n}/4}}$ \\ \hline
$\alpha_{12}=\frac{16\pi^2}{y^2_{\ell n}}$ & 
$u_{12}=\frac{\sin(4\pi y/y_{\ell n})}{\sqrt{nx_{\ell n}y_{\ell n}/2}}$ &
$\alpha_{12}=\frac{36\pi^2}{n^2x^2_{\ell n}}+\frac{4\pi^2}{y^2_{\ell n}}$ & 
$u_{12}=\frac{\sin(\frac{6\pi x}{nx_{\ell n}}+\frac{2\pi y}{y_{\ell n}})}
{\sqrt{nx_{\ell n}y_{\ell n}/4}}$ \\ \hline
$\alpha_{13}=\frac{16\pi^2}{y^2_{\ell n}}$ & 
$u_{13}=\frac{\cos(4\pi y/y_{\ell n})}{\sqrt{nx_{\ell n}y_{\ell n}/2}}$ &
$\alpha_{13}=\frac{36\pi^2}{n^2x^2_{\ell n}}+\frac{4\pi^2}{y^2_{\ell n}}$ & 
$u_{13}=\frac{\cos(\frac{6\pi x}{nx_{\ell n}}+\frac{2\pi y}{y_{\ell n}})}
{\sqrt{nx_{\ell n}y_{\ell n}/4}}$ \\ \hline
$\alpha_{14}=\frac{36\pi^2}{n^2x^2_{\ell n}}$ & 
$u_{14}=\frac{\sin(6\pi x/nx_{\ell n})}{\sqrt{nx_{\ell n}y_{\ell n}/2}}$ &
$\alpha_{14}=\frac{36\pi^2}{n^2x^2_{\ell n}}+\frac{4\pi^2}{y^2_{\ell n}}$ & 
$u_{14}=\frac{\sin(\frac{6\pi x}{nx_{\ell n}}+\frac{2\pi y}{y_{\ell n}})}
{\sqrt{nx_{\ell n}y_{\ell n}/4}}$ \\ \hline
$\alpha_{15}=\frac{36\pi^2}{n^2x^2_{\ell n}}$ & 
$u_{15}=\frac{\cos(6\pi x/nx_{\ell n})}{\sqrt{nx_{\ell n}y_{\ell n}/2}}$ &
$\alpha_{15}=\frac{36\pi^2}{n^2x^2_{\ell n}}+\frac{4\pi^2}{y^2_{\ell n}}$ & 
$u_{15}=\frac{\cos(\frac{6\pi x}{nx_{\ell n}}+\frac{2\pi y}{y_{\ell n}})}
{\sqrt{nx_{\ell n}y_{\ell n}/4}}$ \\ \hline
$\alpha_{16}=\frac{16\pi^2}{n^2x^2_{\ell n}}+\frac{4\pi^2}{y^2_{\ell n}}$ & 
$u_{16}=\frac{\sin(\frac{4\pi x}{nx_{\ell n}}+\frac{2\pi y}{y_{\ell n}})}
{\sqrt{nx_{\ell n}y_{\ell n}/2}}$ &
$\alpha_{16}=\frac{16\pi^2}{n^2x^2_{\ell
n}}+\frac{16\pi^2}{y^2_{\ell n}}$ & 
$u_{16}=\frac{\sin(\frac{4\pi x}{nx_{\ell n}}+\frac{4\pi y}{y_{\ell n}})}
{\sqrt{nx_{\ell n}y_{\ell n}/4}}$ \\ \hline 
$\alpha_{17}=\frac{16\pi^2}{n^2x^2_{\ell n}}+\frac{4\pi^2}{y^2_{\ell n}}$ & 
$u_{17}=\frac{\cos(\frac{4\pi x}{nx_{\ell n}}+\frac{2\pi y}{y_{\ell n}})}
{\sqrt{nx_{\ell n}y_{\ell n}/2}}$ &
$\alpha_{17}=\frac{16\pi^2}{n^2x^2_{\ell
n}}+\frac{16\pi^2}{y^2_{\ell n}}$ & 
$u_{17}=\frac{\cos(\frac{4\pi x}{nx_{\ell n}}+\frac{4\pi y}{y_{\ell n}})}
{\sqrt{nx_{\ell n}y_{\ell n}/4}}$   
\ed{tabular}}{The first 17 eigenvalues and eigenfunctions of $-\Delta
u_i=\alpha_iu_i$.}
\ed{tabela}

The first of the following two lemmas follows from the variational 
characterization for eigenvalues, and the second follows from Lemma 3, the 
Courant nodal domain theorem, 
and geometric properties of the surfaces $W_{\ell/n}$:  

\bg{lema} {\em (\cite{R1}, \cite{R2})}.  Choose $\mu\in{\Bbb Z}^+$ so that 
$\alpha_{\rho_{\ell/n}(\mu)}<4H$. Then {\em Ind}$(W_{\ell/n})\geq\mu-1$.  
\ed{lema} 


\bg{lema} {\em (\cite{R1}, \cite{R2})}.  For all 
$n\in{\Bbb Z}^+$, $n\geq 2$ we 
have that {\em Ind}$(W_{\ell/n})\geq 2n-2$ if $\ell$ is odd, and 
{\em Ind}$(W_{\ell/n})\geq n-2$ if $\ell$ is even.  
\ed{lema}

%
%

\section{The lower bound $8$ for Ind($W_{\ell/n}$)} 

We now show the following: 
\begin{center}
{\bf Numerical Result:}
$\mbox{Ind}(W_{\ell/n})\geq 8$ for all $\ell/n$.
\end{center}

Observe that, although the eigenvalues of $\cal L$ depend on the choice of 
$H$, the number of negative 
eigenvalues is independent of $H$.  So without loss of generality we 
fix $H=1/2$.

By Lemma 5, Ind($W_{\ell/n}$) can be less than $8$ only if $\ell/n$ is one 
of ${3/2}$, ${4/3}$, ${5/3}$, ${5/4}$, ${7/4}$, ${6/5}$, ${8/5}$, 
${8/7}$, ${10/7}$, ${12/7}$, ${10/9}$, $14/9$, or $16/9$.  
Lemma 4 also gives explicit lower bounds for the index, since we know the 
values of $x_{\ell n}$ and $y_{\ell n}$ numerically by formula 
(5), and hence we 
know the $\alpha_{\rho_{\ell/n}(i)}$ (see Table 1).  Lemma 4 implies that the 
index is at least $8$ when $\ell/n$ is $5/4$, $6/5$, $8/7$, $10/7$, or 
$10/9$.  Thus we only need to consider the following eight surfaces: 
\begin{center}
$W_{3/2}$, $W_{4/3}$, $W_{5/3}$, $W_{7/4}$, $W_{8/5}$, $W_{12/7}$, 
$W_{14/9}$, and $W_{16/9}$.
\end{center}
For these surfaces we list, in Table 2, the corresponding $\theta,\
x_{\ell n},\ y_{\ell n}$ and lower bounds for index. These approximate values 
for $\theta$, $x_{\ell n}$, and $y_{\ell n}$ were computed numerically 
using formulas (4) and (5) and the software Mathematica.  
Recall that always $\overline{\theta} \cong 65.354^\circ$.  

\bg{tabela}{
\bg{tabular}{|c|c|c|c|c|c|} \hline
$W_{\ell/n}$ & $\theta$ & $x_{\ell n}$ & $y_{\ell n}$ &
\multicolumn{1}{|p{2.3cm}|}{Lemma 4\newline 
lower bound\newline 
for $\mbox{Ind}(W_{\ell n})$} &
\multicolumn{1}{|p{2.3cm}|}{Lemma 5\newline 
lower bound\newline 
for $\mbox{Ind}(W_{\ell n})$} 
\\ \hline
$W_{3/2}$ & 17.7324$^\circ$ & 2.5556 & 4.2131 & 2 & 2 \\ \hline
$W_{4/3}$ & 12.7898$^\circ$ & 3.2767 & 6.3355 & 6 & 1 \\ \hline
$W_{5/3}$ & 21.4807$^\circ$ & 1.7557 & 2.6402 & 2 & 4 \\ \hline
$W_{7/4}$ & 22.8449$^\circ$ & 1.3315 & 1.9447 & 2 & 6 \\ \hline
$W_{8/5}$ & 20.1374$^\circ$ & 2.0842 & 3.2321 & 2 & 3 \\ \hline
$W_{12/7}$ & 22.3044$^\circ$ & 1.5150 & 2.2380 & 2 & 5 \\ \hline
$W_{14/9}$ & 19.1243$^\circ$ & 2.2970 & 3.6514 & 4 & 7 \\ \hline
$W_{16/9}$ & 23.2182$^\circ$ & 1.1872 & 1.7208 & 2 & 7 \\ \hline
\ed{tabular}}{$x_{\ell n},\ y_{\ell n}$ are computed using the value
$H=1/2$.}
\ed{tabela}

We will find specific spaces on which ${\cal L}$ is
negative definite, for these eight surfaces.  

Let $N$ be an arbitrary positive integer.  Consider a finite subset
$\{\tilde{u}_1=u_{i_1},\dots,\tilde{u}_N=u_{i_N}\}$ of the eigenfunctions
$u_i$ of $-\Delta$ on ${\Bbb C}/\Gamma$, defined in Section 3, with
corresponding eigenvalues $\tilde{\alpha}_j=\alpha_{i_j},\ j=1,\dots,N$. If
we consider any
$u=\sum^N_{i=1}a_i\tilde{u}_i \in\mbox{span}\{\tilde{u}_1,\dots,
\tilde{u}_N\},\ a_1,\dots,a_N \in {\Bbb R}$, then
$\int_{{\Bbb C}/\Gamma}u{\cal L}(u)dx\,dy=
\sum^N_{i,j=1}a_i(\tilde{\alpha}_j\delta_{ij}-\tilde{b}_{ij})a_j$, 
where $\tilde{b}_{ij}:=\int_{{\Bbb C}/\Gamma}V\tilde{u}_i\tilde{u}_j\,dxdy$. 
So we have $\int_{{\Bbb C}/\Gamma}u{\cal L}(u)dx\,dy<0$ for all nonzero
$u\in\mbox{span}\{\tilde{u}_1,\dots,\tilde{u}_N\}$ if and only 
if the matrix
$(\tilde{\alpha}_j\delta_{ij}-\tilde{b}_{ij})_{i,j=1,\dots,N}$
is negative definite.  Lemma 3 then implies:

\bg{teo} {\em (\cite{R1}, \cite{R2})}.  If the $N\times N$ matrix
$(\tilde{\alpha}_j\delta_{ij}-\tilde{b}_{ij})_{i,j=1,\dots,N}$ is negative
definite, then {\em Ind}$(W_{\ell/n})\geq N-1$.
\ed{teo}

\bg{tabela}{
\bg{tabular}{|c|c|c|c|c|c|c|c|c|c|} \hline
$W_{\ell/n}$ & $\tilde{u}_1$ & $\tilde{u}_2$ & $\tilde{u}_3$ &
$\tilde{u}_4$ & $\tilde{u}_5$ & $\tilde{u}_6$ & $\tilde{u}_7$ &
$\tilde{u}_8$ & $\tilde{u}_9$  \\ \hline
$W_{3/2}$ & $u_1$ & $u_2$ & $u_3$ & $u_4$ & $u_5$ & $u_7$ & $u_8$ &
$u_9$ & $u_{17}$  \\ \hline
$W_{4/3}$ & $u_1$ & $u_2$ & $u_3$ & $u_4$ & $u_5$ & $u_6$ & $u_7$ &
$u_8$ & $u_9$  \\ \hline 
$W_{5/3}$ & $u_1$ & $u_2$ & $u_3$ & $u_5$ & $u_6$ & $u_7$ & $u_8$ &
$u_9$ & $u_{15}$  \\ \hline 
$W_{7/4}$ & $u_1$ & $u_2$ & $u_3$ & $u_6$ & $u_7$ & $u_8$ & $u_9$ &
$u_{14}$ & $u_{15}$  \\ \hline 
$W_{8/5}$ & $u_1$ & $u_2$ & $u_3$ & $u_4$ & $u_5$ & $u_{10}$ & $u_{11}$ &
$u_{12}$ & $u_{13}$  \\ \hline 
$W_{12/7}$ & $u_1$ & $u_2$ & $u_3$ & $u_4$ & $u_5$ & $u_{10}$ & $u_{11}$ &
$u_{12}$ & $u_{13}$  \\ \hline 
$W_{14/9}$ & $u_1$ & $u_2$ & $u_3$ & $u_4$ & $u_5$ & $u_{10}$ & $u_{11}$ &
$u_{12}$ & $u_{13}$  \\ \hline 
$W_{16/9}$ & $u_1$ & $u_2$ & $u_3$ & $u_4$ & $u_5$ & $u_{10}$ & $u_{11}$ &
$u_{12}$ & $u_{13}$ \\ \hline    
\ed{tabular}}{Eigenfunctions of $-\Delta$ producing $9$-dimensional 
spaces on which $Q$ is negative definite.} 
\ed{tabela}

\bg{dft} 
{\em Given $A,B$ even integers and $\ell,n \in {\Bbb Z}^+$, we now define 
the following basic integrals: 
\begin{center}
$I_0(\ell,n,A,B)=\frac{1}{nx_{\ell n}y_{\ell n}}
\int^{x_{\ell n}/4}_0\int^{y_{\ell n}/4}_0V
\left(\cos\frac{2\pi x}{x_{\ell n}}\right)^A
\left(\cos\frac{2\pi y}{y_{\ell n}}\right)^B dydx$, 
\end{center}\begin{center}
$I_1(\ell,n)=\frac{8}{nx_{\ell n}y_{\ell n}}
\int^{nx_{\ell n}/4}_0\int^{y_{\ell n}/4}_0V
\cos\left(\frac{4\pi x}{nx_{\ell n}}\right) dydx$,
\end{center}\begin{center}
$I_2(\ell,n)=\frac{8}{nx_{\ell n}y_{\ell n}}
\int^{nx_{\ell n}/4}_0\int^{y_{\ell n}/4}_0V
\cos\left(\frac{8\pi x}{nx_{\ell n}}\right) dydx$,
\end{center}\begin{center}
$I_3(\ell,n)=\frac{8}{nx_{\ell n}y_{\ell n}}
\int^{nx_{\ell n}/4}_0\int^{y_{\ell n}/4}_0V
\cos\left(\frac{16\pi x}{nx_{\ell n}y_{\ell n}}\right) dydx$,
\end{center}\begin{center}
$I_4(\ell,n)=\frac{8}{nx_{\ell n}y_{\ell n}}
\int^{nx_{\ell n}/4}_0\int^{y_{\ell n}/4}_0V
\cos\left(\frac{4\pi x}{nx_{\ell n}}\right)
\cos\left(\frac{4\pi y}{y_{\ell n}}\right) dydx$,
\end{center}\begin{center}
$I_5(\ell,n)=\frac{8}{nx_{\ell n}y_{\ell n}}
\int^{nx_{\ell n}/4}_0\int^{y_{\ell n}/4}_0V
\cos\left(\frac{4\pi x}{nx_{\ell n}}\right)
\cos\left(\frac{8\pi x}{nx_{\ell n}}\right)dydx$,
\end{center}\begin{center}
$I_6(\ell,n)=\frac{8}{nx_{\ell n}y_{\ell n}}
\int^{nx_{\ell n}/4}_0\int^{y_{\ell n}/4}_0V
\cos\left(\frac{8\pi x}{nx_{\ell n}}\right)
\cos\left(\frac{4\pi y}{y_{\ell n}}\right)dydx$,
\end{center}\begin{center}
$I_7(\ell,n)=\frac{8}{nx_{\ell n}y_{\ell n}}
\int^{nx_{\ell n}/4}_0\int^{y_{\ell n}/4}_0V
\cos\left(\frac{12\pi x}{nx_{\ell n}}\right)
\cos\left(\frac{4\pi y}{y_{\ell n}}\right)dydx$. 
\end{center}
}
\ed{dft}

Now, for each surface $W_{\ell/n}$ given in Table 2, we will fix $N=9$
and choose the subset $\{\tilde{u}_1,\dots,\tilde{u}_9\}$ such that
the matrix $(\tilde{\alpha}_j\delta_{ij}-\tilde{b}_{ij})_{i,j=1,\dots,N}$ is
negative definite. These choices are given in Table 3.  
With these choices for $\tilde{u}_i$, we have the following lemma:  

\bg{lema} With the choices given in Table 3, all elements of the 
eight matrices 
${\cal M}(\ell,n):=(\tilde{\alpha}_j\delta_{ij}-
\tilde{b}_{ij})_{i,j=1,\dots,9}$ can be 
expressed in terms of the basic integrals 
$I_0(\ell,n,A,B)$ and $I_j(\ell,n)$ for 
$A,B$ even and $j=1,2,\dots,7$.

\medskip

\bgpf{} {\em The symmetries 
$V(x,y)=V(-x,y)=V(x,-y)=V\left(\frac{x_{\ell
n}}{2}-x,y\right)=V\left(x,\frac{y_{\ell n}}{2}-y\right)$ of $V$ and 
the identities $\cos(a\pm b)=\cos(a)\cos(b)\mp\sin(a)\sin(b),\ \sin(a\pm 
b)=\sin(a)\cos(b)\pm\sin(b)\cos(a)$ give the relations shown in Table 4, 
proving the lemma.}
\ed{lema}

\bg{tabela}{
\bg{tabular}{|c|p{14cm}|} \hline
\multicolumn{1}{|p{1cm}|}{For\newline 
$W_{3/2}$} & $M_{1,1}=-32\ I_0(3,2,0,0)$,
$M_{4,4}=\alpha_4-64(I_0(3,2,0,0)-I_0(3,2,0,2))$\newline
$M_{5,5}=\alpha_5-64\ I_0(3,2,0,2)$,
$M_{7,7}=\alpha_7-64\ I_0(3,2,2,0)$\newline
$M_{17,17}=\alpha_{17}-64(I_0(3,2,0,0)-I_0(3,2,0,2)
                 -I_0(3,2,2,0)+2\ I_0(3,2,2,2))$\newline
$M_{i,i}=\alpha_i-32I_0(3,2,0,0)$, for $i=2,3,8,9$\newline
$M_{i,j}=0$ for $i < j$ with $i,j \in \{ 1,2,3,4,5,7,8,9,17 \}$.\\ \hline
\multicolumn{1}{|p{1cm}|}{For\newline 
$W_{4/3}$} & $M_{1,1}=-48\ I_0(4,3,0,0)$, 
$M_{2,2}=\alpha_2-48\ I_0(4,3,0,0)+2\ I_2(4,3)$\newline
$M_{3,3}=\alpha_3-48\ I_0(4,3,0,0)-2\ I_2(4,3)$, 
$M_{4,4}=\alpha_4-48\ I_0(4,3,0,0)+2\ I_4(4,3)$\newline
$M_{5,5}=\alpha_5-48\ I_0(4,3,0,0)-2\ I_4(4,3)$, 
$M_{6,6}=\alpha_6-48\ I_0(4,3,0,0)+2\ I_4(4,3)$\newline
$M_{7,7}=\alpha_7-48\ I_0(4,3,0,0)-2\ I_4(4,3)$\newline
$M_{8,8}=\alpha_8-384\ (I_0(4,3,0,2)-\ I_0(4,3,0,4))$\newline
$M_{9,9}=\alpha_9-96\ (4\ I_0(4,3,0,4)-4\ I_0(4,3,0,2)+I_0(4,3,0,0))$\newline
$M_{1,3}=-2\sqrt{2}I_1(4,3)$, 
$M_{1,9}=-48\sqrt{2}(-I_0(4,3,0,0)+2I_0(4,3,0,2))$\newline
$M_{4,6}=-48(-I_0(4,3,0,0)+2I_0(4,3,0,2))+2I_1(4,3)$\newline
$M_{5,7}=-48(-I_0(4,3,0,0)+2I_0(4,3,0,2))-2I_1(4,3)$\newline
$M_{3,9}=-4\ I_4(4,3)$, 
all other $M_{i,j}$ with $i < j \leq 9$ are zero.\\ \hline
\multicolumn{1}{|p{1cm}|}{For \newline 
$W_{5/3}$} & $M_{1,1}=-48\ I_0(5,3,0,0)$, 
$M_{2,2}=\alpha_2-48\ I_0(5,3,0,0)+2\ I_1(5,3)$\newline
$M_{3,3}=\alpha_3-48\ I_0(5,3,0,0)-2\ I_1(5,3)$, 
$M_{5,5}=\alpha_5-96\ I_0(5,3,0,2)$\newline
$M_{6,6}=\alpha_6-48\ I_0(5,3,0,0)+2\ I_2(5,3)$, 
$M_{7,7}=\alpha_7-48\ I_0(5,3,0,0)-2\ I_2(5,3)$\newline
$M_{8,8}=\alpha_8-48\ I_0(5,3,0,0)+2\ I_4(5,3)$, 
$M_{9,9}=\alpha_9-48\ I_0(5,3,0,0)-2\ I_4(5,3)$\newline
$M_{15,15}=\alpha_{15}-96\ I_0(5,3,2,0)$, 
$M_{1,7}=-2\sqrt{2}I_1(5,3)$\newline
$M_{3,15}=-2(I_1(5,3)+I_2(5,3))$, 
all other pertinent $M_{i,j}$ with $i < j$ are zero.\\ \hline
\multicolumn{1}{|p{1cm}|}{For\newline 
$W_{7/4}$} & $M_{i,i}=\alpha_i-64\ I_0(7,4,0,0)$ 
for $i=1,2,3,6,7,8,9,14,15$\newline
all other pertinent $M_{i,j}$ with $i < j$ are zero.\\ \hline
\multicolumn{1}{|p{1cm}|}{For\newline 
$W_{8/5}$,\newline $W_{12/7}$,\newline $W_{14/9}$,\newline $W_{16/9}$} & 
$M_{1,1}=-16n\ I_0(\ell,n,0,0)$, 
$M_{2,2}=\alpha_2-16n\ I_0(\ell,n,0,0)+2I_1(\ell,n)$\newline
$M_{3,3}=\alpha_3-16n\ I_0(\ell,n,0,0)-2I_1(\ell,n)$, 
$M_{4,4}=\alpha_4-16n\ I_0(\ell,n,0,0)+2I_4(\ell,n)$\newline
$M_{5,5}=\alpha_5-16n\ I_0(\ell,n,0,0)-2I_4(\ell,n)$\newline
$M_{10,10}=\alpha_{10}-16n\ I_0(\ell,n,0,0)+2I_3(\ell,n)$\newline
$M_{11,11}=\alpha_{11}-16n\ I_0(\ell,n,0,0)-2I_3(\ell,n)$\newline
$M_{12,12}=\alpha_{12}-16n\ I_0(\ell,n,0,0)+2I_7(\ell,n)$\newline
$M_{13,13}=\alpha_{13}-16n\ I_0(\ell,n,0,0)-2I_7(\ell,n)$, 
$M_{1,3}=-2\sqrt{2}I_1(\ell,n)$\newline
$M_{1,11}=-2\sqrt{2}I_2(\ell,n)$, 
$M_{2,10}=-4\ I_1(\ell,n)+4\ I_5(\ell,n)$\newline
$M_{3,11}=-4\ I_5(\ell,n)$, 
$M_{4,12}=-2I_1(\ell,n)+2I_6(\ell,n)$\newline
$M_{5,13}=-2I_1(\ell,n)-2I_6(\ell,n)$, 
all other pertinent $M_{i,j}$ with $i < j$ are zero.\\ \hline
\ed{tabular}}
{Elements $M_{i,j}$ of the symmetric matrices 
${\cal M}(\ell,n)$ expressed in terms of the basic integrals.  We have 
chosen here to index the $M_{i,j}$ using the counters associated to 
$\alpha_j$ and $u_j$, rather than $\tilde{\alpha}_j$ and $\tilde{u}_j$.}
\ed{tabela}

By numerical methods, we can estimate that all of the relevant 
$I_j(\ell,n)$ for $j \geq 1$ are approximately zero, and that 
\begin{center}
$ I_0(3,2,0,0)\cong0.2968$, $I_0(3,2,2,0)\cong0.2304$, 
$I_0(3,2,0,2)\cong0.2408$, $I_0(3,2,2,2)\cong0.1947$, 
\end{center}\begin{center} 
$I_0(4,3,0,0)\cong0.1077$, $I_0(4,3,0,2)\cong0.0776$, 
$I_0(4,3,0,4)\cong0.0667$, $I_0(5,3,0,0)\cong0.4532$, 
\end{center}\begin{center} 
$I_0(5,3,2,0)\cong0.3910$, $I_0(5,3,0,2)\cong0.4046$, 
$I_0(7,4,0,0)\cong0.6072$, $I_0(8,5,0,0)\cong0.1878$, 
\end{center}\begin{center} 
$I_0(12,7,0,0)\cong0.2652$, $I_0(14,9,0,0)\cong0.0841$, 
$I_0(16,9,0,0)\cong0.3419$. 
\end{center}
These values were computed with a Mathematica program using the 
NIntegrate and JacobiCN commands, and the program is available at the 
web site of the third author.  One note of warning is that Mathematica 
has different conventions than Walter's paper, and hence $cn_k$ in 
\cite{Wa} is equivalent to $cn_{k^2}$ in Mathematica.  We include a sample 
of our code in the Appendix.  

Now we can make approximations for the eight matrices 
${\cal M}(\ell,n)$.  

The matrix ${\cal M}(3,2)$ is approximately 
\[ {\cal M}(3,2) \approx \left(
\begin{array}{ccccccccc}
-9.50&0&0&0&0&0&0&0&0 \vspace{-0.05in} \\ 
0&-7.99&0&0&0&0&0&0&0 \vspace{-0.05in} \\
0&0&-7.99&0&0&0&0&0&0 \vspace{-0.05in} \\
0&0&0&-1.36&0&0&0&0&0 \vspace{-0.05in} \\
0&0&0&0&-13.2&0&0&0&0 \vspace{-0.05in} \\
0&0&0&0&0&-8.70&0&0&0 \vspace{-0.05in} \\
0&0&0&0&0&0&-5.76&0&0 \vspace{-0.05in} \\
0&0&0&0&0&0&0&-5.76&0 \vspace{-0.05in} \\
0&0&0&0&0&0&0&0&-5.50
\end{array}
\right) \; , \]
and all nondiagonal terms are known to be zero by rigorous mathematical 
computation, and all nonzero entries have been computed only numerically.  

${\cal M}(4,3)$ is approximately the nondiagonal matrix 
\[ {\cal M}(4,3) \hspace{-0.014in} \approx \hspace{-0.02in} \left(
\begin{array}{cccccccccc}
-5.17&0&{\cal O}&0&0&0&0&0&-3.23 \vspace{-0.05in} \\
0&-3.53&0&0&0&0&0&0&0  \vspace{-0.05in} \\
{\cal O}&0&-3.53&0&0&0&0&0&{\cal O}  \vspace{-0.05in} \\
0&0&0&-3.78&0&-2.29&0&0&0  \vspace{-0.05in} \\
0&0&0&0&-3.78&0&-2.29&0&0  \vspace{-0.05in} \\
0&0&0&-2.29&0&-3.78&0&0&0  \vspace{-0.05in} \\
0&0&0&0&-2.29&0&-3.78&0&0  \vspace{-0.05in} \\
0&0&0&0&0&0&0&-0.25&0  \vspace{-0.05in} \\
-3.23&0&{\cal O}&0&0&0&0&0&-2.21
\end{array}\right) \; , \] and again here all entries that are $0$ have 
been computed mathematically rigorously, and all nonzero entries have been 
computed only numerically.  The symbol $\cal O$ denotes an entry that has 
been computed numerically to be approximately 
zero, but not mathematically rigorously.
We shall continue to use these conventions in all remaining matrices.  

${\cal M}(5,3)$ is approximately 
the diagonal matrix 
\[ {\cal M}(5,3) \hspace{-0.014in} \approx \hspace{-0.02in} \left(
\begin{array}{cccccccccc}
-21.8&0&0&0&0&{\cal O}&0&0&0 \vspace{-0.05in} \\
0&-20.3&0&0&0&0&0&0&0  \vspace{-0.05in} \\
0&0&-20.3&0&0&0&0&0&{\cal O}  \vspace{-0.05in} \\
0&0&0&-33.2&0&0&0&0&0  \vspace{-0.05in} \\
0&0&0&0&-16.1&0&0&0&0  \vspace{-0.05in} \\
{\cal O}&0&0&0&0&-16.1&0&0&0  \vspace{-0.05in} \\
0&0&0&0&0&0&-14.7&0&0  \vspace{-0.05in} \\
0&0&0&0&0&0&0&-14.7&0  \vspace{-0.05in} \\
0&0&{\cal O}&0&0&0&0&0&-24.7
\end{array}\right) \; . \] 

${\cal M}(7,4)$ is approximately 
the diagonal matrix 
\[ {\cal M}(7,4) \hspace{-0.014in} \approx \hspace{-0.02in} \left(
\begin{array}{cccccccccc}
-38.9&0&0&0&0&0&0&0&0 \vspace{-0.05in} \\
0&-37.5&0&0&0&0&0&0&0  \vspace{-0.05in} \\
0&0&-37.5&0&0&0&0&0&0  \vspace{-0.05in} \\
0&0&0&-33.3&0&0&0&0&0  \vspace{-0.05in} \\
0&0&0&0&-33.3&0&0&0&0  \vspace{-0.05in} \\
0&0&0&0&0&-27.0&0&0&0  \vspace{-0.05in} \\
0&0&0&0&0&0&-27.0&0&0  \vspace{-0.05in} \\
0&0&0&0&0&0&0&-26.3&0  \vspace{-0.05in} \\
0&0&0&0&0&0&0&0&-26.3
\end{array}\right) \; . \] 

${\cal M}(8,5)$ is approximately 
the diagonal matrix 
\[ {\cal M}(8,5) \hspace{-0.014in} \approx \hspace{-0.02in} \left(
\begin{array}{cccccccccc}
-15.0&0&{\cal O}&0&0&0&{\cal O}&0&0 \vspace{-0.05in} \\
0&-13.6&0&0&0&{\cal O}&0&0&0  \vspace{-0.05in} \\
{\cal O}&0&-13.6&0&0&0&{\cal O}&0&0  \vspace{-0.05in} \\
0&0&0&-10.9&0&0&0&{\cal O}&0  \vspace{-0.05in} \\
0&0&0&0&-10.9&0&0&0&{\cal O}  \vspace{-0.05in} \\
0&{\cal O}&0&0&0&-9.2&0&0&0  \vspace{-0.05in} \\
{\cal O}&0&{\cal O}&0&0&0&-9.2&0&0  \vspace{-0.05in} \\
0&0&0&{\cal O}&0&0&0&-8.0&0  \vspace{-0.05in} \\
0&0&0&0&{\cal O}&0&0&0&-8.0
\end{array}\right) \; . \] 

${\cal M}(12,7)$ is approximately 
the diagonal matrix 
\[ {\cal M}(12,7) \hspace{-0.014in} \approx \hspace{-0.02in} \left(
\begin{array}{cccccccccc}
-29.7&0&{\cal O}&0&0&0&{\cal O}&0&0 \vspace{-0.05in} \\
0&-28.3&0&0&0&{\cal O}&0&0&0  \vspace{-0.05in} \\
{\cal O}&0&-28.3&0&0&0&{\cal O}&0&0  \vspace{-0.05in} \\
0&0&0&-21.5&0&0&0&{\cal O}&0  \vspace{-0.05in} \\
0&0&0&0&-21.5&0&0&0&{\cal O}  \vspace{-0.05in} \\
0&{\cal O}&0&0&0&-24.1&0&0&0  \vspace{-0.05in} \\
{\cal O}&0&{\cal O}&0&0&0&-24.1&0&0  \vspace{-0.05in} \\
0&0&0&{\cal O}&0&0&0&-18.7&0  \vspace{-0.05in} \\
0&0&0&0&{\cal O}&0&0&0&-18.7
\end{array}\right) \; . \] 

${\cal M}(14,9)$ is approximately 
the diagonal matrix 
\[ {\cal M}(14,9) \hspace{-0.014in} \approx \hspace{-0.02in} \left(
\begin{array}{cccccccccc}
-12.1&0&{\cal O}&0&0&0&{\cal O}&0&0 \vspace{-0.05in} \\
0&-11.7&0&0&0&{\cal O}&0&0&0  \vspace{-0.05in} \\
{\cal O}&0&-11.7&0&0&0&{\cal O}&0&0  \vspace{-0.05in} \\
0&0&0&-9.1&0&0&0&{\cal O}&0  \vspace{-0.05in} \\
0&0&0&0&-9.1&0&0&0&{\cal O}  \vspace{-0.05in} \\
0&{\cal O}&0&0&0&-10.6&0&0&0  \vspace{-0.05in} \\
{\cal O}&0&{\cal O}&0&0&0&-10.6&0&0  \vspace{-0.05in} \\
0&0&0&{\cal O}&0&0&0&-8.3&0  \vspace{-0.05in} \\
0&0&0&0&{\cal O}&0&0&0&-8.3
\end{array}\right) \; . \] 

${\cal M}(16,9)$ is approximately 
the diagonal matrix 
\[ {\cal M}(16,9) \hspace{-0.014in} \approx \hspace{-0.02in} \left(
\begin{array}{cccccccccc}
-49.2&0&{\cal O}&0&0&0&{\cal O}&0&0 \vspace{-0.05in} \\
0&-47.9&0&0&0&{\cal O}&0&0&0  \vspace{-0.05in} \\
{\cal O}&0&-47.9&0&0&0&{\cal O}&0&0  \vspace{-0.05in} \\
0&0&0&-35.6&0&0&0&{\cal O}&0  \vspace{-0.05in} \\
0&0&0&0&-35.6&0&0&0&{\cal O}  \vspace{-0.05in} \\
0&{\cal O}&0&0&0&-43.7&0&0&0  \vspace{-0.05in} \\
{\cal O}&0&{\cal O}&0&0&0&-43.7&0&0  \vspace{-0.05in} \\
0&0&0&{\cal O}&0&0&0&-32.8&0  \vspace{-0.05in} \\
0&0&0&0&{\cal O}&0&0&0&-32.8
\end{array}\right) \; . \] 

All eight of these matrices are $9 \times 9$ and negative definite. Hence
Theorem 1 implies the numerical result.

\section{Appendix: the Mathematica code}

The following is a Mathematica code for computing the values 
$I_0(4,3,0,0)$, $I_0(4,3,0,2)$, $I_0(4,3,0,4)$, $I_1(4,3)$, 
$I_2(4,3)$, $I_4(4,3)$, and the elements of the matrix ${\cal M}(4,3)$.
The seven other needed codes for different $\ell$ and $n$ were written 
similarly.  

{\small 
\begin{verbatim}
H = 1/2;     k1 = Sin[theta1];     k2 = Sin[theta2];
gamma1 = Sqrt[Tan[theta1]];     gamma2 = Sqrt[Tan[theta2]];
alpha1 = Sqrt[4 H Sin[2theta2]/Sin[2(theta1 + theta2)]];
alpha2 = Sqrt[4 H Sin[2theta1]/Sin[2(theta1 + theta2)]];
F = 4ArcTanh[gamma1 gamma2 JacobiCN[alpha1 x, k1^2] JacobiCN[alpha2 y, k2^2]];
V = 4 H Cosh[F];

ell = 4;     n = 3;
theta1 = 2 Pi (12.7898/360);     theta2 = 2 Pi (65.354955354/360); 
x0 = 3.2767;     y0 = 6.3355;

Print["I_0(4,3,0,0) is ",I0x4c3c0c0x = (1/(n x0 y0)) NIntegrate[
   V , {x, 0, x0/4}, {y, 0, y0/4}]];

Print["I_0(4,3,0,2) is ",I0x4c3c0c2x = (1/(n x0 y0)) NIntegrate[
   V (Cos[2 Pi x/x0])^0(Cos[2 Pi y/y0])^2,{x, 0, x0/4}, {y, 0, y0/4}]];

Print["I_0(4,3,0,4) is ",I0x4c3c0c4x = (1/(n x0 y0)) NIntegrate[
   V (Cos[2 Pi x/x0])^0(Cos[2 Pi y/y0])^4,{x, 0, x0/4}, {y, 0, y0/4}]];      

Print["I_1(4,3) is ",I1x4c3x = (8/(n x0 y0)) (NIntegrate[
   V (Cos[4 Pi x/(n x0)]),{x, 0, x0/4}, {y, 0, y0/4}] + NIntegrate[
   V (Cos[4 Pi x/(n x0)]),{x, x0/4, 2 x0/4}, {y, 0, y0/4}] + NIntegrate[
   V (Cos[4 Pi x/(n x0)]),{x, 2 x0/4, n x0/4}, {y, 0, y0/4}])];

Print["I_2(4,3) is ",I2x4c3x = (8/(n x0 y0)) (NIntegrate[
   V (Cos[8 Pi x/(n x0)]),{x, 0, x0/4}, {y, 0, y0/4}] + NIntegrate[
   V (Cos[8 Pi x/(n x0)]),{x, x0/4, 2 x0/4}, {y, 0, y0/4}] + NIntegrate[
   V (Cos[8 Pi x/(n x0)]),{x, 2 x0/4, n x0/4}, {y, 0, y0/4}])];

Print["I_4(4,3) is ",I4x4c3x = (8/(n x0 y0)) (NIntegrate[
  V (Cos[4 Pi x/(n x0)]) (Cos[4 Pi y/y0]),{x,0,x0/4},{y,0,y0/4}]+NIntegrate[
  V (Cos[4 Pi x/(n x0)]) (Cos[4 Pi y/y0]),{x,x0/4,2 x0/4},{y,0,y0/4}]+NIntegrate[
  V (Cos[4 Pi x/(n x0)]) (Cos[4 Pi y/y0]),{x, 2 x0/4, n x0/4}, {y, 0, y0/4}])];

aa = 0;  bb = 0;     alpha1 = aa(4 N[Pi^2]/(n^2 x0^2)) + bb (4 N[Pi^2]/(y0^2));

aa = 4;  bb = 0;     alpha2 = aa(4 N[Pi^2]/(n^2 x0^2)) + bb (4 N[Pi^2]/(y0^2));
alpha3 = alpha2;

aa = 1;  bb = 1;     alpha4 = aa(4 N[Pi^2]/(n^2 x0^2)) + bb (4 N[Pi^2]/(y0^2));
alpha5 = alpha4;     alpha6 = alpha4;     alpha7 = alpha4;

aa = 0;  bb = 4;     alpha8 = aa(4 N[Pi^2]/(n^2 x0^2)) + bb (4 N[Pi^2]/(y0^2));
alpha9 = alpha8;

Print["M(1,1) is ", -48 I0x4c3c0c0x];
Print["M(2,2) is ", alpha2 - 48 I0x4c3c0c0x];
Print["M(3,3) is ", alpha3 - 48 I0x4c3c0c0x];
Print["M(4,4) is ", alpha4 - 48 I0x4c3c0c0x];
Print["M(5,5) is ", alpha5 - 48 I0x4c3c0c0x];
Print["M(6,6) is ", alpha6 - 48 I0x4c3c0c0x];
Print["M(7,7) is ", alpha7 - 48 I0x4c3c0c0x];
Print["M(8,8) is ", alpha8 - 384 (I0x4c3c0c2x - I0x4c3c0c4x)];
Print["M(9,9) is ", alpha9 - 96 (I0x4c3c0c0x + 4 I0x4c3c0c4x - 4 I0x4c3c0c2x)];
Print["M(1,9) is ", -48 N[Sqrt[2]] (-I0x4c3c0c0x + 2 I0x4c3c0c2x)];
Print["M(4,6) is ", -48 (-I0x4c3c0c0x + 2 I0x4c3c0c2x)];
Print["M(5,7) is ", -48 (-I0x4c3c0c0x + 2 I0x4c3c0c2x)];
\end{verbatim}
}

{\small 
\bg{thebibliography}{999}

{\small 

\bibitem[A]{A} U. ABRESCH, {\it Constant mean curvature tori in terms
of elliptic functions}, J. reine u. angew Math. {\bf374} (1987), 169-192.

\bibitem[BC]{BC} L. BARBOSA, M. do CARMO, {\it Stability of
hypersurfaces with constant mean curvature}, Math. Z. {\bf185} (1984), 
339--353.


\bibitem[CP]{CP} M. do CARMO, C. K. PENG, {\it Stable minimal surfaces in 
${\Bbb R}^3$ are planes}, Bull. Amer. Math. Soc. {\bf1} (1979), 
903--906.  



\bibitem[FC]{FC} D. FISCHER-COLBRIE, {\it On complete minimal surfaces 
with finite Morse index in three manifolds}, Invent. Math. {\bf82} (1985), 
121--132.

\bibitem[H]{H} H. HOPF, {\it Differential geometry in the large}, 
Lecture Notes in Math. {\bf1000} Springer-Verlag (1983).  

\bibitem[LR]{LR} F. J. LOPEZ, A. ROS, {\it Completenimal surfaces with 
index one and stable constant mean curvature surfaces}, Comment. Math. 
Helvetici {\bf64} (1989), 34--43.  

\bibitem[R1]{R1} W. ROSSMAN, {\it The Morse index of Wente tori}, to appear 
in Geom. Dedicata.

\bibitem[R2]{R2} W. ROSSMAN, {\it Wente Tori and Morse Index}, An. Acad. 
Bras. Ci. {\bf 71(4)} (1999), 607--613.  

\bibitem[Si]{Si} A. M. SILVEIRA, {\it Stability of complete noncompact 
surfaces with constant mean curvature}, Math. Ann. {\bf277} (1987), 
629--638.  

\bibitem[Sp]{Sp} J. SPRUCK, {\it The elliptic sinh-Gordon equation and the 
construction of toroidal soap bubbles}, 
Lecture Notes in Math. {\bf1340} Springer-Verlag (1988), 275--301.  

\bibitem[Wa]{Wa} R. WALTER, {\it Explicit examples to the $H$-problem 
of Heinz Hopf}, Geometriae Dedicata {\bf23} (1987), 187--213.

\bibitem[We]{We} H. WENTE, {\it Counterexample of a conjecture of H.
Hopf}, Pacific J. Math. {\bf121} (1986), 193--243.

}

\ed{thebibliography}
}
\vskip1cm

\parbox[t]{6in}{
Levi Lopes de Lima, Departamento de Matem\'{a}tica, Universidade
Federal do Cear\'{a} \\Campus do Pici, 60455--760 Fortaleza,
Brazil, levi@mat.ufc.br}

\bigskip

\parbox[t]{6in}{
Vicente Francisco de Sousa Neto, Departamento de Matem\'atica\\ 
Universidade Cat\'olica de Pernambuco, Recife, Brazil, 
vicente@unicap.br}

\bigskip

\parbox[t]{6in}{
Wayne Rossman, Department of Mathematics, Faculty of Science\\ Kobe
University, Rokko, Kobe 657-8501,
Japan, wayne@math.kobe-u.ac.jp\\
http://www.math.kobe-u.ac.jp/HOME/wayne/wayne.html}

\ed{document}